\documentclass[10pt]{article}

\input xy
\xyoption{all}
\usepackage{amsmath}
\usepackage{amsfonts}
\usepackage{amssymb}
\usepackage{amsthm,stmaryrd}
\usepackage{graphicx,epsfig}
\usepackage{bbm,color}
\usepackage{dbnsymb}
\newtheorem{lem}{Lemma}
\newtheorem{quest}{Question}
\newtheorem{prop}{Proposition}
\newtheorem{rem}{Remark}

\newcommand{\boA}{\mathcal{A}}
\newcommand{\boB}{\mathcal{B}}
\newcommand{\boD}{\mathcal{D}}
\newcommand{\boH}{\mathcal{H}}
\newcommand{\boM}{\mathcal{M}}
\newcommand{\Q}{\mathbb{Q}}
\newcommand{\alg}[3]{\vphantom{#2}_{#1} #2_{#3}}

\newcommand{\arete}{\!\!\frown}
\newcommand{\Zr}{Z^{\rm rat}}
\newcommand{\Zwr}{Z^{\fourwheel {\rm rat}}}
\DeclareMathOperator{\lift}{Lift}
\DeclareMathOperator{\hair}{Hair}
\begin{document}
\title{On Kontsevich Integral of torus knots
\footnote{\emph{keywords}: finite type invariants, Kontsevich integral, torus knots, Wheels and Wheeling, rationality. \quad \quad \quad \quad \quad \quad \quad \emph{2000 Mathematics Subject Classification:} 57M27.}}
\author{Julien March\'e
\footnote{Institut de Math\'ematiques de Jussieu, \'Equipe ``Topologie et G\'eom\'etries Alg\'ebriques''
Case 7012, Universit\'e Paris VII, 75251 Paris CEDEX 05, France. 
\quad e-mail: \texttt{marche@math.jussieu.fr}}}
\date{October 2003}
\maketitle

\begin{abstract}

We study the unwheeled rational Kontsevich integral of torus knots. We give a precise formula for these invariants up to loop degree 3 and show that they appear as colorings of simple diagrams. We show that they behave under cyclic branched coverings in a very simple way. Our proof is combinatorial: it uses the results of Wheels and Wheelings and new decorations of diagrams.

\end{abstract}

\section{Introduction and notations}

This article is divided in 4 parts: the first one explains the notations used in the sequel and some facts known about rationality, the second part is a computation of the unwheeled Kontsevich integral of torus knots up to loop degree 3 using formal series which encode some series of diagrams. Then, in the third part, we compute a rational form of the preceding expression and show that it appears as a coloring of chain diagrams as it is suggested by the figure \ref{loop}. As a consequence of this computation, we show that the operator $\lift_r$ which corresponds to cyclic branched coverings of $S^3$ along the knot simply acts on a diagram $D$ of loop degree lower than 3 by multiplying it by $r^{-\chi(D)}$ where $\chi$ is the Euler characteristic.

\begin{figure}[htbp]
\begin{center}
\input{loop.pstex_t}
\caption{Diagrams appearing in the unwheeled Kontsevich integral of 
torus knots}
\label{loop}
\end{center}
\end{figure}

The initial idea for this computation is not new, it has been used by Christine Lescop (see \cite{les}) and Dror Bar-Natan in an unpublished work. Lev Rozansky also computed formulas for the loop expansion of torus knots in the weight system associated to $\rm{sl}_2$ (see \cite{roz2}).
The computation of the 2-loop part of torus knots has been done independantly by Tomotada Ohtsuki in \cite{oht} who computed more generally a formula for 2-loop part of knots cabled by torus knots.
We would like to thank Stavros Garoufalidis for useful remarks and Marcos Marino for pointing out a mistake in the last formula of this article.

\subsection{Normalizations of the Kontsevich integral}

Let $K$ be a knot in $S^3$ and suppose that $K$ has a banded structure with self-linking 0. We will note $Z(K)$ the Kontsevich integral of $K$ in the algebra $\boA$ of trivalent diagrams lying on a circle.

Let $\boB$ be the algebra of uni-trivalent diagrams. It is well known that the Poincar\'e-Birkhoff-Witt map $\chi:\boB\to\boA$ is an isomorphism but not an algebra isomorphism. We will note $\sigma$ its inverse.

If $U$ is the trivial knot, we define $\Omega=\sigma Z(U)$. The map $\Upsilon=\chi \circ \partial_{\Omega}:\boB\to \boA$ defined by instance in \cite{th} is known to be an algebra isomorphism. The quantity $Z^{\fourwheel}(K)=\Upsilon^{-1}Z(K)$ will be called unwheeled Kontsevich integral and behaves better than $\sigma Z(K)$ under connected sum and cyclic branched coverings.

For each knot $K$, the quantities $Z(K)$, $\sigma Z(K)$ and $Z^{\fourwheel}(K)$ are group-like, which means that they are exponentials of a series of connected diagrams. We will note respectively $z(K)$, $\sigma z(K)$ and $z^{\fourwheel}(K)$ the logarithm of these quantities.

\subsection{Loop degree and rationality}

If $D$ is a connected diagram of $\boB$, its first Betti number defines a degree called loop degree. The loop degree 1 part of $\sigma Z(K)$ or $Z^{\fourwheel}(K)$ is well-known: it only depends on the Alexander polynomial of $K$. For the higher degrees, very few is known. There are formulas for the 2-loop part of small knots in Rozansky's table (see \cite{roz}), and we can find in \cite{double} a formula for the 2-loop part of untwisted whitehead doubles. In the sequel, we give a formula for the 2-loop and 3-loop parts of the Kontsevich integral of torus knots.

In order to make precise computations, we will need the following formalism:
let $\boH$ be a cocommutative Hopf algebra (or a Hopf algebra up to completion) and $\boM$ be an algebra over $\boH$. We consider a space of diagrams noted  $\boD(\boH,\boM)$ which was defined in \cite{vog}. It is roughly obtained by decorating the edges of a diagram by elements of $\boM$ and allowing elements of $H$ to slide through vertices thanks to the coproduct law.

We know that the space $\boB$ is isomorphic to $\boD(\Q[[h]],\Q[[h]])$. In particular, the diagrams of loop degree 1 appear as colorings of the circle by an even power series without constant term. In the following, we will call wheels such diagrams and identify a power series with the wheel series it represents.

Let $f(x)$ be the power series defined by $\frac{1}{2}\log{\frac{\sinh{x/2}}{x/2}}$. The famous wheel formula (see \cite{th}) states that $\sigma z(U)=f(x)$. Further, it was shown in \cite{kr} that the loop degree 1 part of $\sigma z(K)$ is $f(x)+Wh_K(x)$ where $Wh_K(x)=-\frac{1}{2}\log \Delta(e^x)$ and $\Delta$ is the Alexander polynomial of $K$.

As we are interested in the higher loop degree part, we need to recall the rationality theorem which was proved in \cite{rat}:

Let $\Lambda=\Q[t,t^{-1}]$ and $\Lambda_{\rm loc}$ the localization of $\Lambda$ with respect to elements $f$ satisfying $f(1)=1$.
The substitution $t=\exp(h)$ gives a morphism between couples  $(\Lambda,\Lambda_{\rm loc})$ and $(\Q[[h]],\Q[[h]])$, hence a morphism between $\boD(\Lambda,\Lambda_{\rm loc})$ and $\boB$. 

We call $\hair$ this application: the rationality theorem tells us that the series  $\sigma z(K)$ minus the loop degree 1 part lies in the image of the $\hair$ map. More precisely, there is an element $\Zr(K)\in \boD(\Lambda,\Lambda_{\rm loc})$ whose denominators on each edge is at most $\Delta(K)$ such that 
$\sigma Z(K)=\exp(f(x)+Wh_K(x)) \hair \Zr(K)$. Respectively, there is such an element for the unwheeled invariant but we will give the formula later.

A construction, developed in \cite{rat} gives $\Zr(K)$ as an invariant of $K$ which is not automatic because the $\hair$ map is not injective, although it is in small degrees (see \cite{ber}).

\section{Computation of the torus knot integral up to loop degree 3}

Let $p$ and $q$ be two coprime integers such that $p>0$. We note $K_{p,q}$ the torus banded knot with parameters $p$ and $q$ and self-linking 0, and $L_{p,q}$ the torus banded knot with banding parallel to the torus on which it lies. This knot has self-linking $pq$, and his Kontsevich integral is a bit easier to compute. 

The method of computation is inspired from \cite{les}: we first compute the Kontsevich integral of the following braid.

Let $p$ points be lying on the vertices of a regular $p$-gone. We note $\gamma$ be the braid obtained by rotating the whole picture by an angle $2\pi\frac{q}{p}$.

Let associate to any one dimensional manifold $\Gamma$ the space $\boA(\Gamma)$ of trivalent diagrams lying on $\Gamma$. It defines a contravariant functor with respect to continuous maps relative to boundaries. Let $\phi_p^*$ be the map induced by the projection on the first factor $\phi_p:[0,1]\times\{1,\ldots,p\}\to [0,1]$ and $\isolatedchord$ be the only degree 1 diagram in $\boA([0,1])$. 
Then a direct computation of monodromy of the K-Z connection shows that $\gamma$ has a Kontsevich integral equal to $\phi_p^*( \exp_{\#}(\frac{q}{2p}\isolatedchord))$.

The banded knot $L_{p,q}$ is obtained by closing the previous braid: this translates diagrammatically to the following: let $\psi_p$ be the map from $S^1$ to itself defined by $\psi_p(z)=z^p$. Then, $Z(L_{p,q})=\psi_p^*(\nu\# \exp_{\#}(\frac{q}{2p}\isolatedchord))$, where $\nu=Z(U)$.

By lemma 4.10 of \cite{th}, the map $\psi_p^*$ viewed in $\boB$ has the following form:
if $D\in \boB$ has $k$ legs (i.e. univalent vertices) then $\sigma \psi_p^* \chi D  = p^k D$. We will note more simply $D_p$ the result of this operation which looks like a change of variable.

Then, to compute $Z(K_{p,q})$ from $Z(L_{p,q})$, we only need to change the framing, that is $Z(K_{p,q})=\exp_{\#}(-\frac{pq}{2}\isolatedchord)\#Z(L_{p,q})$. We will transform this product in the usual one by applying the unwheeling map $\Upsilon^{-1}$. As a result, we will have a formula for $Z^{\fourwheel}(K_{p,q})$.

We now sum up the steps of the computation:

\begin{enumerate}
\item Computation of $\sigma(\nu\# \exp_{\#}(\frac{q}{2p}\isolatedchord))$
\item Change of variables $x\mapsto px$
\item Unwheeling
\end{enumerate}

To state the result of our computations, we will need a way to present some diagrams of loop degree lower than 3. The first diagram will code 1-loop part, the second, 2-loop part and the two last diagrams 3-loop part. 

\begin{itemize}
\item We recall that $x^n$ is a wheel with $n$ legs ($n$ is even).

\item Let $x^n y^m$ represent two wheels glued on one edge, with $n$ remaining legs on the left and $n$ remaining legs on the right. In particular, $n$ and $m$ are odd and $x^n y^m=x^m y^n$.

\item The expression ${z_1}^{m_1}x^n {z_2}^{m_2}$ represents the coloring of a diagram with three wheels joined by two edges, with $n$ legs on the central wheel and $m_1$ and $m_2$ legs on the other wheels. By convention, we put the variable associated to the middle wheel between the two others.

\item
Finally, we note $[x^{n},x^{m}]$ the sum of all diagrams obtained by gluing two wheels of size $n$ and $m$ in two points.
\end{itemize}

\begin{figure}[htbp]
\begin{center}
\input{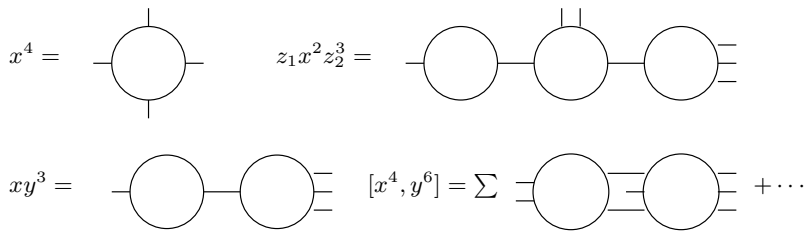}
\caption{Examples of diagrams}
\label{roues}
\end{center}
\end{figure}

The remaining of this section is devoted to the proof of the following proposition:

\begin{prop}
Up to loop degree 3, the unwheeled Kontsevich integral of torus knot can be expressed by the following power series, where $c(x)=f(px)+f(qx)-f(pqx)$.
\begin{eqnarray*}
c(x) & \text{1-loop part}&\\
f'(py)f'(qx)-c'(x)f'(pqy)& \text{2-loop part}&\\
\frac{1}{2}f'(pz_1)f''(qx)f'(pz_2)+\frac{1}{2}f'(q z_1)f''(px)f'(q z_2)&&\quad\quad\quad\quad\quad\quad(*)\\
+\frac{1}{2}f'(pqz_1)c''(x)f'(pqz_2)-\frac{1}{2}c'(z_1)f''(pqx)c'(z_2)\\
-f'(pq z_1)[p f''(px)f'(qz_2)+q f''(qx)f'(pz_2)-pqf''(pqz_2)c'(x)]\\
\frac{1}{p^2q^2}[f(qx),f(px)]-\frac{1}{p^2q^2}[f(pqx),c(x)]
\end{eqnarray*}
\end{prop}

\begin{rem}
From this formula, we see that the 2-loop part appears as a coloring of a dumb-bell graph. But often, such a series of diagram is presented as a coloring of the Theta graph. Indeed, given a coloring of the dumb-bell $x^ny^m$, we can make an (IHX) move on the central edge, and replace this diagram by two times a graph $\ThetaGraph$ with $n$ legs on the upper edge, and $m$ legs on the lower edge but in the opposite direction. Not all colorings of $\ThetaGraph$ are obtained in that way but we can easily express our formulas as colorings of $\ThetaGraph$.
\end{rem}

\subsection{Step 1:}

From this point, we will be interested only in diagrams with loop degree 1, 2 or 3 as we were not able till now to perform a computation for general loop degree. 

It is clear that all the operations considered below as $\#$, $\partial$, etc... cannot decrease the loop degree. Then it is licit to quotient the spaces of diagrams by all diagrams with loop degree greater than 3.

There is a less licit simplification we make: we also quotient by all diagrams without legs or not lying on a one dimensional manifold. Indeed, the unwheeled invariants may have such elements and we will justify this choice in the section \ref{closed}. 

We recall that $\Upsilon=\chi\circ\partial_{\Omega}$ is an algebra isomorphism and that $\Upsilon^{-1}\nu = \Omega$ and $\Upsilon^{-1}\isolatedchord= \frown$ (modulo closed diagrams). 

Then, what we have to compute is
 $$\sigma(\nu\# \exp_{\#}(\frac{q}{2p}\isolatedchord))=
\partial_{\Omega} (\Omega\exp(\frac{q}{2p}\arete)).$$

As this expression is group-like, we just have to find all connected diagrams appearing in it.

In the computation of $\partial_{\Omega} (\Omega\exp(\frac{q}{2p}\arete))$, we will concentrate on the wheels of the derived term and make the following observation: such a wheel can only be glued to another wheel, and each gluing increase the loop degree by one. Hence, there can be at most two gluings. 
We want to separate the diagrams obtained by gluing these wheels in two points. To compute these contributions, we state the following very useful lemma.

\begin{lem}
\label{lemme1}
\begin{itemize}
\item
Let $f(x)$ and $g(y)$ be two power series of wheels. The series obtained by gluing them in one point is $f'(x)g'(y)$.
\item
Let $f(x)$ and $g(y)h(z)$ be two series of diagrams. The series obtained by gluing them in one point is $f'(x)g'(y)h(z)+f'(x)h'(z)g(y)$.
\end{itemize}
\end{lem}

The proof is straightforward and has a first application in the following lemma:

\begin{lem}
\label{lemme2}
The expression obtained by gluing the derived wheels in two points is
$$\frac{p^2}{q^2}[f(\frac{q}{p}x),f(x)]+\frac{1}{2}f'(\frac{q}{p}z_1)f''(x)f'(\frac{q}{p}z_2).$$
\end{lem}

\begin{proof}
Recall that $\Omega=\exp(f(x))$. The first term comes from the gluing of a right wheel to the same left wheel. We have to fill the remaining left legs with the $\frac{q}{2p}$ struts in the only way such that it does not increase the loop degree. Finally, we just multiply by $\frac{q}{p}$ to the power the number of remaining legs. There are as many legs as the left wheel, minus two, which justifies the first part of the formula.

For the second term, we just apply twice the lemma \ref{lemme1}.
\end{proof}

We can suppose now that each wheel can be glued in at most one leg. Let us examine precisely the following expression:

$$\partial_{\Omega} (\Omega\exp(\frac{q}{2p}\arete))=\langle \Omega_x,\Omega_{x+y}\exp(\frac{q}{2p}\alg{x+y}{\frown}{x+y})\rangle_x.$$ 

We can replace $\exp(\omega_{x+y})$ by the series $\exp(\omega_y +\hat{\omega}_x)$. The series $\hat{\omega}_x$ is obtained by replacing a wheel by the sum of the coloring of its legs by one $x$ (and the other ones $y$). 

We use a famous trick concerning the doubling of the strut part to get:

\begin{eqnarray*}
\partial_{\Omega}(\exp(\omega+\frac{q}{2p}\arete))
&=
\langle \Omega_x, \exp( \omega_y +\hat{\omega_x}+
\frac{q}{2p}\alg{x}{\arete}{x}+\frac{q}{p}\alg{x}{\arete}{y}+\frac{q}{2p}\alg{y}{\arete}{y}) \rangle_x \\
&=
\langle \partial_{\exp(\hat{\omega}_x)}\partial_{\exp(\frac{q}{2p}\alg{x}{\arete}{x})}  \Omega_x,\exp(\frac{q}{p}\alg{x}{\arete}{y})\rangle_x \exp(\frac{q}{2p}\alg{y}{\arete}{y})\Omega_y\rangle_x \\
&=
\langle \partial_{\exp(\hat{\omega}_x)} \Omega_x,\exp(\frac{q}{p}\alg{x}{\arete}{y})\rangle_x \exp(\frac{q}{2p}\alg{y}{\arete}{y})\Omega_y
\end{eqnarray*}

But as $\hat{\omega}_{x+x'}=\hat{\omega}_{x}+\hat{\omega}_{x'}$, the operator $\partial_{\exp(\hat{\omega}_x)}$ commutes with exponentials and then 
$$\partial_{\exp(\hat{\omega}_x)} \Omega_x=\exp(\partial_{\exp(\hat{\omega}_x)}\omega_x)=\exp(f(x)+f'(y)f'(x)+\frac{1}{2}f'(z_1)f''(x)f'(z_2)).$$

The final connected contribution is $\frac{q}{2p}\arete+ f(x)+ f(\frac{q}{p}x)+f'(y)f'(\frac{q}{p}x)+\frac{1}{2}f'(z_1)f''(\frac{q}{p}x)f'(z_2)$.

Putting all terms together, we find the following expression for $\partial_{\Omega} (\Omega\exp(\frac{q}{2p}\arete))$:

$$\exp(\frac{q}{2p}\arete+f(x)+f(\frac{q}{p}x)+f'(y)f'(\frac{q}{p}x)+\frac{1}{2}f'(z_1)f''(\frac{q}{p}x)f'(z_2)+
\frac{1}{2}f'(\frac{q}{p}z_1)f''(x)f'(\frac{q}{p}z_2)+
\frac{p^2}{q^2}[f(\frac{q}{p}x),f(x)])$$

\subsection{Step 2:}
This step is very simple. In order to compute $\sigma Z(L_{p,q})$, we just have to multiply each term by as many factors $p$ as legs.

\begin{align*}
\sigma Z(L_{p,q})=&
\exp(\frac{pq}{2}\arete+f(px)+f(qx)+f'(py)f'(qx) \\
&+\frac{1}{2}f'(pz_1)f''(qx)f'(pz_2)+
\frac{1}{2}f'(q z_1)f''(px)(q z_2)+
\frac{1}{p^2q^2}[f(qx),f(px)])
\end{align*}

\subsection{Step 3:}

In this section, we want to correct the framing defect. The only way we know is to unwheel the preceding expression.
Let us note $c(x)=f(px)+f(qx)-f(pqx)$ and factorize it. We then compute

$$Z^{\fourwheel}(K_{p,q})=\partial_{\Omega^{-1}}\exp(\frac{pq}{2}\arete+f(pqx)+c(x)+R)\exp(-\frac{pq}{2}\arete).$$
Here, $R=f'(py)f'(qx)+f'(pz_1)f''(qx)f'(pz_2)+
\frac{1}{2}f'(q z_1)f''(px)f'(q z_2)+
\frac{1}{p^2q^2}[f(qx),f(px)]$ is a series of loop degree 2 and 3.

As before, we first compute the contribution coming from the gluing of $c$ in two points of the same wheel. The lemma \ref{lemme2} give us the following expression, noted $U$:

$$U=-\frac{1}{p^2q^2}[f(pqx),c(x)]+\frac{1}{2}f'(pqz_1)c''(x)f'(pqz_2).$$

If we suppose that $c$ can be glued to only one leg, we see that the following diagram is commutative:

$$
\xymatrix{
\boB \ar[rr]^{\times\exp(-c)}\ar[dd]^{\partial_{\Omega}} &&\boB\ar[dd]^{\partial_{\tilde{\Omega}}} \\ \\
\boB \ar[rr]^{\times\exp(-c)}&&\boB
}
$$
Where $\tilde{\Omega}=\partial_{\exp(\hat{c})}\Omega=\exp(f(x)+c'(x)f'(y)+\frac{1}{2}c'(z_1)f''(x)c'(z_2))$.

But 
$$\partial_{\tilde{\Omega}}\exp(\frac{pq}{2}\frown)=\langle \partial_{\exp(\hat{c}_x)}\Omega_x,\exp(\frac{pq}{2}\alg{x+y}{\frown}{x+y})\rangle_x={\tilde{\Omega}}_{pq}\exp(\frac{pq}{2}\frown).$$

If we factorize this element in the expression of $Z^{\fourwheel}(K_{p,q})$ we get:

$$\partial_{\tilde{\Omega}^{-1}}\left[ \partial_{\tilde{\Omega}}\exp(\frac{pq}{2}\frown) \exp(R-c'(x)f'(pqy)-\frac{1}{2}c'(z_1)f''(pqx)c'(z_2))\right]\exp(-\frac{pq}{2}\arete+c+U).$$

Now, we just have to take care of the action by derivation of $\tilde{\Omega}^{-1}$ on the 2-loop part of the right member, that is:

$$-f'(pq z_1)[p f''(px)f'(qz_2)+q f''(qx)f'(pz_2)-pqf''(pqz_2)c'(x)].$$

Finally, we collect all terms in increasing loop-degree order and conclude the computation.

\section{Rationality}

In this part, we will express the diagrams of formula $(*)$ in a rational form and show the following proposition:
\begin{prop}
Let $X_n$ be the connected diagram consisting of $n$ wheels connected by $n-1$ edges. Then for $n\le 3$, the unwheeled rational invariant of $K_{p,q}$ of loop degree $n$ appear as a coloring of $X_n$.
\end{prop}

\subsection{Loop degree 1 and 2}

It is easy to check that the 1-loop part is as expected:

\begin{align*}
c(x)-f(x)&=f(px)+f(qx)-f(pqx)-f(x)\\
&=\frac{1}{2}\log\frac{\sinh(px/2)}{px/2}+\frac{1}{2}\log\frac{\sinh(qx/2)}{qx/2}-\frac{1}{2}\log\frac{\sinh(pqx/2)}{pqx/2}-\frac{1}{2}\log\frac{\sinh(x/2)}{x/2}\\
&=\frac{1}{2}\log\frac{\sinh(px/2)\sinh(qx/2)}{\sinh(pqx/2)\sinh(x/2)} =-\frac{1}{2}\log D_{p,q}(e^x)
\end{align*}
Here, $D_{p,q}(t)=t^{-\frac{1}{2}(p-1)(q-1)}\frac{(t^{pq}-1)(t-1)}{(t^p-1)(t^q-1)}$ is the Alexander polynomial of $K_{p,q}$.

In the following, we will need a small extension of the spaces of diagrams, because it happens that we need to use rational expressions with poles at the unity. Let us note $\Lambda'$ the field of fractions of $\Lambda=\Q[t,t^{-1}]$. This is an algebra over $\Lambda$ and we have an injective $\Lambda$-morphism from $\Lambda_{\rm loc}$ to $\Lambda'$.

Then we have a map from $\boD(\Lambda,\Lambda_{\rm loc})\to\boD(\Lambda,\Lambda')$. There is a corresponding $\hair$ map which fits in the following diagram:
$$
\xymatrix{
\boD(\Lambda,\Lambda_{\rm loc}) \ar[r] \ar[d]^{\hair}& \boD(\Lambda,\Lambda')\ar[d]^{\hair}\\
\boD(\Q[[h]],\Q[[h]]) \ar[r] & \boD(\Q[[h]],\Q[[h]][h^{-1}])
}$$

Although none of the previous maps is injective, we will identify all diagrams with their image in the latter space (these maps happen to be injective in low degrees).

We can now compute the 2-loop term which we call $z_2$, using the formula $f'(x)=\frac{1}{4}\coth(x/2)-\frac{1}{2x}$, all terms containing $\frac{1}{x}$ simplify after symetrization. If we note $t=e^x$ and $r=e^y$, then
\begin{eqnarray*}
z_2=\frac{1}{32}[2pq\frac{t^{pq}+1}{t^{pq}-1}\frac{s^{pq}+1}{s^{pq}-1}-p\frac{t^{pq}+1}{t^{pq}-1}\frac{s^{p}+1}{s^{p}-1}-q\frac{t^{pq}+1}{t^{pq}-1}\frac{s^{q}+1}{s^{q}-1}\\
-p\frac{s^{pq}+1}{s^{pq}-1}\frac{t^{p}+1}{t^{p}-1}-q\frac{s^{pq}+1}{s^{pq}-1}\frac{t^{q}+1}{t^{q}-1}+\frac{t^{p}+1}{t^{p}-1}\frac{s^{q}+1}{s^{q}-1}+\frac{t^{q}+1}{t^{q}-1}\frac{s^{p}+1}{s^{p}-1}]
\end{eqnarray*}

This formula appears as a coloring of the dumb-bell graph with denominators dividing $t^{pq}-1$. If we want to write it with denominators dividing $D_{p,q}$, we are forced to write $z_2$ as a coloring of $\ThetaGraph$. We can show that there is such a factorization but we were not able to find a close formula for the numerators.

\subsection{Loop degree 3}

Let us study the 3-loop term, noted $z_3$. We decompose it in two parts, $z_3^1$ and $z_3^2$. 
The $z_3^2$ part is just the part of $z_3$ expressed with brackets (the last line of the formula $(*)$).

In the expression of $z_3^1$, if we write $f''(x)=-\frac{1}{8}\frac{1}{\sinh(x/2)^2}+\frac{1}{2x^2}$ and develop, we get a sum of two terms. The first does not contain any fractional term and is obtained from $z_3^1$ by forgetting all of them (and hence is rational). Concerning the second term, a computation with MAPLE shows that it reduces to: 

$$z'_3=-\frac{1}{2pq}\left[  f'(py)\frac{1}{x^2}f'(qz)-c'(y)\frac{1}{x^2}f'(pqz)\right].$$

\begin{prop}
The series $z'3$ and $z_3^2$ cancels.
\end{prop}

\begin{proof}

Let us compute the term $z_3^2$:

We recall that $z_3^2=\frac{1}{p^2q^2}[f(qx),f(px)]-\frac{1}{p^2q^2}[f(pqx),c(x)]$. 
It is obtained from $z_2$ by summing all gluings of a left leg on a right leg. There is a normalization factor $\frac{1}{pq}$ and $\frac{1}{2}$ which counterbalance the order given to the two gluings.

We will need the following two lemma which interpret some diagrams with inverse legs.

\begin{lem}
\label{lem1}
Let $D$ be a diagram consisting of a circle glued on a segment by an edge. Imagine that the circle is colored by a series $g(x)$, and the segment by $\frac{1}{x}$. Using IHX relations we can make the moves suggested by the figure \ref{laurent}. Making the series sliding, we can cancel the $\frac{1}{x}$ term except for one term which is just the opposite of the initial term. 

This shows that the initial diagram can be expressed by a coloring of another diagram without inverse legs.

\begin{figure}[htbp]
\begin{center}
\input{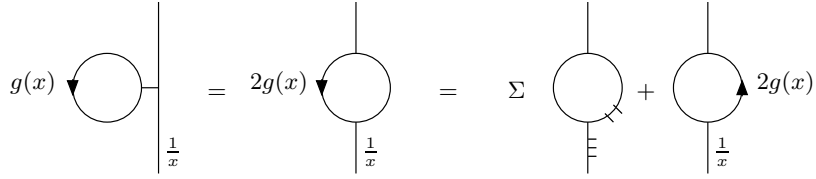}
\caption{Diagram containing a $\frac{1}{x}$ term}
\label{laurent}
\end{center}
\end{figure}

\end{lem}

\begin{lem}
\label{lem2}
Let $g$ be a series coloring a circle attached to an edge.

Consider the sum of the diagrams obtained by gluing the end of a free edge to the legs defined by $g$.
Then this series is obtained by the diagram of figure \ref{laurent2}.

\begin{figure}[htbp]
\begin{center}
\input{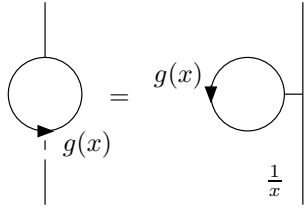}
\caption{Gluing a leg and a wheel}
\label{laurent2}
\end{center}
\end{figure}
\end{lem}

\begin{proof}
We just have to check it for $g(x)=x^n$. 
In the left hand side figure, note $y$ a leg lying on the right part of the circle and $z$ a leg on the vertical segment ($z=x-y$).

The left hand side figure is just obtained by the coloring $\sum_{i+j=n-1}x^i y^j$.
Concerning the right part, as in the lemma \ref{lem1}, we make $x^n$ slide to get $(z+y)^n$, eliminate the term $y^n$, and divide by $z$. We compute $[(z+y)^n-y^n]/z=(x^n-y^n)/(x-y)=\sum_{i+j=n-1}x^i y^j$. This proves the lemma.
\end{proof}

If we apply the lemma \ref{lem2} to both sides of the 2-loop part $z_2$, we show that $z_3^2$ and $z'_3$ cancel.
\end{proof}

We saw in this computation that both 2-loop terms and 3-loop terms appear as colorings of very special diagrams. The following question is quite natural:

\begin{quest}
Which kind of diagram are needed in higher loop degrees to express the rational Kontsevich invariant of torus knot? It is unlikely that the degree $n$ part is a coloring of $X_n$, but we can certainly reduce the number of diagrams needed.
\end{quest}

\subsection{The case of closed diagrams}\label{closed}

We recall that in all these computations we neglected closed diagrams, although they exist in the expression of $Z^{\fourwheel}(K_{p,q})$. 
We propose to show that we do not need to add any closed diagram to the expression of $\Zwr(K_{p,q})=\exp(c+z_2+z_3)$.

First, there are no closed diagram in the 1-loop or 2-loop terms. The only one comes from $z_3$, and precisely
from the diagram colored by $z_1\frac{1}{x^2}z_2$. A direct computation show that this term is 
$\frac{(p^2-1)(q^2-1)}{1152}\ThetaGraph_2$. In the sequel, $\ThetaGraph_n$ is the graph $\ThetaGraph$ with the middle edge replaced by n parallel copies.

Following \cite{lift}, the invariant $Z^{\fourwheel rat}$ is normalized such that for loop degree greater than 1,
$$\hair^{\frac{1}{\langle \Omega,\Omega \rangle}}\Zwr(K)=Z^{\fourwheel}(K)$$
As $\sigma Z(K_{p,q})$ does not contain any closed diagram, we know that $\langle Z^{\fourwheel}(K), \Omega \rangle=1$. We must then show that $\langle \hair \Zwr(K_{p,q}), \Omega \rangle = \langle \Omega, \Omega \rangle $.

But in the left hand side, the only possible gluing comes from the 1-loop part $\frac{1}{48}(p^2+q^2-p^2q^2)\twowheel$ of $c$ and $\frac{1}{48}\twowheel$ of $\Omega$. Summed with the closed diagram coming from $z_3$, we get $\frac{1}{1152}\ThetaGraph_2=\langle \Omega,\Omega\rangle$.

\section{Branched coverings}

\subsection{Main formula}

A great interest for rational expression of Kontsevich integral comes from its relation with branched coverings. More precisely, if $K_{p,q}$ is the torus knot of parameters $p$ and $q$, and $r$ is an integer, let us note $\Sigma^r(K_{p,q})$ be the pair formed of the cyclic branched covering of $S^3$ of order $r$ over $K_{p,q}$ and the ramification link.

If $r$ is coprime with $p$ and $q$, the ramification locus is a knot, and the underlying 3-manifold is a rational homology sphere, the Brieskorn manifold $\Sigma(p,q,r)$. 

In \cite{lift}, a map $\lift_r$ is described which intertwines rational invariant of the cyclic branched coverings and rational invariant of the initial knot in the following way:

$$Z^{\fourwheel rat}(\Sigma^r(K))= \exp(\frac{\sigma_r(K)}{16}\ThetaGraph)\lift_r Z^{\fourwheel rat}(K).$$

We now study this map in the case of torus knots and loop degree lower than 3 and prove the following proposition:

\begin{prop}
Note $\Zwr(K_{p,q})=\exp(\sum\limits_{k>0} z_k)$ decomposing by loop degree as before,
then for all  $k\le 3$, $$\lift_r z_k = r^{k-1} z_k.$$
\end{prop}

\begin{proof}

For $k=1$, the Alexander polynomial is invariant by these branched coverings. The same will be true for the 1-loop part $z_1$.

Let us look to the other cases.
The $\lift_r$ map was defined only for diagrams decorated by fractions without poles at $r$-roots of unity. As we extended the decorations to all fractions, the definition of $\lift_r$ makes sense for any diagram.
In the definition of the $\lift_r$ map, we need to express all denominators as polynomials of $t^r$. Then, we look to the numerators as a coloring by monomials, which is the same as a linear combination of 1-cohomology classes of the underlying graph. We keep only the classes divisible by $r$ and divide them, then we put back denominators replacing $t^r$ by $t$. Finally we multiply the result by $r$. 
This construction is very easy in our case because every edge colored is part of a circle without other colorings and the colorings are the following:
\begin{itemize}
\item[] let $n$ be coprime with $r$ ($n=p,q$ or $pq$)
\item
$f_1=\frac{t^{n}+1}{t^{n}-1}$ comes from the first derivative $f'(x)=\frac{1}{4}\coth(x/2)-\frac{1}{2x}$
\item
$f_2=
\frac{t^n}{(t^n-1)^2}$ comes from the second derivative $f''(x)=-\frac{1}{8}\frac{1}{\sinh(x/2)^2}+\frac{1}{2x^2}$.
\end{itemize}
In order to express $f_1$ with the right denominator, we multiply numerator and denominator by $1+t^n+\cdots+(t^n)^{r-1}$ to get 
$\frac{1+2t^n+\cdots+2t^{n(r-1)}+t^{nr}}{t^{nr}-1}$. The only $r$-divisible numerators are $1$ and $t^{nr}$ because $n$ and $r$ are coprime. The result of the $\lift_r$ map is then the identity (before multiplying by $r$).

We do the same operation with $f_2$ and multiply numerator and denominator by $(1+t^n+\cdots+(t^n)^{r-1})^2$ to get 
$\frac{t^n(1+t^n+\cdots+(t^n)^{r-1})^2}{(t^{nr}-1)^2}$.
The numerator is $t^n(1+t^n+\cdots+(t^n)^{r-1})^2 =\sum_{i,j=0}^{r-1}t^{n(i+j+1)}$. The indices for which the order of the monomial is $r$-divisible are such that $i+j=n-1$. There are $r$ such terms, and the result of the $\lift_r$ map is then the multiplication by $r$ (before multiplying by $r$ at the end)

For the $z_2$ term, we have two terms $f_1$, and then $\lift_r z_2 = r z_2$.
For the $z_3$ term, we have two terms $f_1$ and one term $f_2$ and then $\lift_r z_3 =r^2 z_3$.
This ends the proof of the proposition.
\end{proof}

\begin{quest}
Can we extend this proposition for larger values of $k$? Do we have any conceptual interpretation of this formula?
\end{quest} 

\subsection{Application to LMO invariant of Brieskorn Spheres}

The LMO invariant of $\Sigma(p,q,r)$ is just the closed part of $\sigma Z(\Sigma^r(K_{p,q}))$ i.e. $\langle Z^{\fourwheel}(\Sigma^r(K_{p,q})),\Omega\rangle$. 
Using the map $\lift$ we have the following formula: $${\rm LMO}(\Sigma(p,q,r))=\exp(\frac{\sigma_r(K)}{16}\ThetaGraph)\frac{\langle \lift_r \Zr(K_{p,q}), \Omega \rangle}{\langle \Omega,\Omega \rangle}.$$

From this formula, we can prove that the degree two term is $\frac{(p^2-1)(q^2-1)(r^2-1)}{1152}\Theta_2$.

Supposing the proposition is true for $n=4$ we deduced  the following formula for the degree 3 term: ${\rm LMO}(\Sigma(p,q,r))_3=\frac{-pqr(p^2-1)(q^2-1)(r^2-1)}{13824}\Theta_3$.

These computations agree with formulas for the LMO invariant of Seiferts spaces which can be found in \cite{bn} and \cite{mm}.

\nocite{*}
\bibliographystyle{halpha}
\bibliography{toric}

\end{document}